\def\PaperDate{2003/10/13}
\documentclass[12pt]{article}
  \usepackage[german,american]{babel}
  \usepackage{datum}
  \usepackage{bux2ref}
  \usepackage{graphicx}
  \usepackage[GlobalNumbered]{buxmath}
  \usepackage[USletter]{buxlayout}

  \usepackage{epic,eepic}
  
  \InputIfFileExists{basic_notation.tex}{}{}
  \InputIfFileExists{notation.tex}{}{}

  \newcommand{\IncludeMathDiagram}[2][]{%
    \raisebox{-0.5\height}{\fbox{\includegraphics[#1]{#2}{}{}}}
  }

  \DeclareFontFamily{OT1}{chess}{}
  \DeclareFontShape{OT1}{chess}{m}{n}{<-> chess30}{}
  \DeclareFixedFont{\chessfont}{OT1}{chess}{m}{n}{70}
  \DeclareTextFontCommand{\chess}{\fontfamily{OT1}{chess}\fontseries{m}\fontshape{n}}

  \renewvariable{\TheDim}{d}
  \renewvariable{\TheProjection}{\pi}
  \newvariable{\TheDummyNumber}{\nu}

  \newcommand{\centergraphics}[2][]{%
    \clap{$\begin{array}{c}\includegraphics[#1]{#2}{}{}\end{array}$}
  }      
\begin{document}
  \title{Braiding and Tangling the Chessboard Complex}
  \author{Kai-Uwe Bux}
  \date{\datum{\PaperDate}}
  \maketitle
  \begin{abstract}
    We describe a series of complexes that relate to the braid groups as the
    matching complexes relate to the symmetric groups. A modified construction
    applies as well to other complexes based on edge sets in graphs. We
    show that our constructions will yield Cohen-Macauley complexes
    provided the underlying complexes are Cohen-Macauley.

    Finally, we discuss a related series of complexes to provide some positive
    evidence that the braided Houghton groups $\HoughtonBraid[\TheTopNumber]$,
    introduced by F.~Degenhardt, are of type \FType[\TheTopNumber-\One] but
    not of type \FType[\TheTopNumber].
  \end{abstract}

  \section{Preliminaries}
    A space $\TheSpace$ is called \notion{$\TheDim$-connected} if
    $\HomotopyOf[\TheIndex]{\TheSpace}=\Zero$ for $\Zero\leq\TheIndex\leq\TheDim$.
    It is called $\Parentheses{-\One}$-connected, if $\TheSpace$ is non-empty.
    A space is
    called \notion{$\TheDim$-spherical} if it of dimension $\TheDim$ and
    $\Parentheses{\TheDim-\One}$-connected. Note that a $\TheDim$-spherical
    complex is homotopy equivalent to a (possibly infinite and possibly empty)
    wedge of $\TheDim$-spheres.
    A CW-complex is \notion{spherical} if it is $\TheDim$-spherical for some
    dimension $\TheDim$.

    The CW-complexes we are concerned with in this paper
    are mostly \notion{$\Delta$-complexes} \cite{Hatcher:2001}.
    These can be viewed as piece-wise Euclidean
    cell complexes whose cells have the shapes
    of regular simplices. The metric information,
    however, is not really essential. You
    might want to think of them just as very nice
    CW-complexes whose cells look like simplices.
    Note that every simplicial complex is a $\Delta$-complex.

    Links in $\Delta$-complexes, in general, do not equal the boundaries of
    stars: A bigon is a perfectly valid $\Delta$-complex, and each of its
    vertices will have a link that consists of precisely two points.
    However, if a $\Delta$-complex happens to be a simplicial complex then
    the two notions of links coincide.

    A poset $\ThePoset$ is called \notion{$\TheDim$-connected}
    or \notion{$\TheDim$-spherical} if its
    geometric realization enjoys said property. For talking about posets, we find a
    topological language convenient: For any element $\ThePosetElement\in\ThePoset$,
    we call
    \begin{notationlist}
      \item
        \(
          \ClosureOf{\ThePosetElement}
          :=
          \SetOf[\AltPosetElement\in\ThePoset]{
            \AltPosetElement \isfaceof \ThePosetElement
          }
        \)
        the \notion{closure} of $\ThePosetElement$,
      \item
        \(
          \BoundaryOf{\ThePosetElement}
          :=
          \SetOf[\AltPosetElement\in\ThePoset]{
            \AltPosetElement \isstrictfaceof \ThePosetElement
          }
        \)
        the \notion{boundary} of $\ThePosetElement$,
      \item
        \(
          \StarOf{\ThePosetElement}
          :=
          \SetOf[\AltPosetElement\in\ThePoset]{
            \ThePosetElement \isfaceof \AltPosetElement
          }
        \)
        the \notion{star} of $\ThePosetElement$, and
      \item
        \(
          \LinkOf{\ThePosetElement}
          :=
          \SetOf[\AltPosetElement\in\ThePoset]{
            \ThePosetElement \isstrictfaceof \AltPosetElement
          }
        \)
        the \notion{link} of $\ThePosetElement$. The intersection
      \item
        \(
          \OpenInterval{\ThePosetElement}{\AltPosetElement}
          :=
          \LinkOf{\ThePosetElement}
          \intersect
          \BoundaryOf{\AltPosetElement}
          =
          \SetOf[\ThirdPosetElement\in\ThePoset]{%
            \ThePosetElement \isstrictfaceof \ThirdPosetElement \isstrictfaceof \AltPosetElement
          }
        \)
        is the \notion{open interval} from $\ThePosetElement$ to
        $\AltPosetElement$.

        The \notion{geometric realization}
      \item
        \(
          \RealizationOf{\ThePoset}
        \)
        is the simplicial complex of finite
        $\isstrictfaceof$-chains in $\ThePoset$.
        The \notion{dimension}
      \item
        \(
          \DimOf{\ThePoset}
        \)
        of a poset $\ThePoset$ is the dimension of its geometric realization.
        If the dimension
        of $\ThePoset$ is finite, a maximum length chain in $\ThePoset$ has
        length $\DimOf{\ThePoset}$.
        The \notion{height}
      \item
        \(
          \HeightOf{\ThePosetElement}
        \)
        of an element $\ThePosetElement$ is the dimension of its closure.
    \end{notationlist}
    This terminology is, of course, inspired by the poset of cells
    in a regular CW-complex
    where $\isstrictfaceof$ is given by the face relation. %
      We mention that people who view posets from a more algebraic
      angle might prefer a different terminology: The closure of an
      element $\ThePosetElement$ is often refered to as the
      \notion{principal order ideal generated by $\ThePosetElement$} and the star
      of $\ThePosetElement$ is often called the \notion{principle
      filter generated by $\ThePosetElement$}. In this
      note, however, a topological terminology seems to be more appropriate.
    
    We follow Quillen's influential paper \cite{Quillen:1978} and call
    a simplicial complex $\TheSimplicialComplex$
    \notion{Cohen-Macauley} if it is spherical and
    every simplex $\TheSimplex$ has a link of dimension
    $\DimOf{\TheSimplicialComplex}-\DimOf{\TheSimplex}-\One$ that
    is spherical, as well. A poset $\ThePoset$
    is \notion{Cohen-Macauley}
    if its \notion{geometric realization} (i.e., the associated simplicial
    complex of chains in $\ThePoset$) is Cohen-Macauley. Quillen
    observes \cite[Proposition~8.6]{Quillen:1978}
    that $\ThePoset$ is Cohen-Macauley if and only if
    it is spherical and all
    links, boundaries, and open intervals in $\ThePoset$ are spherical, too.%
    \footnote{%
      It is necessary to point out that those posets are often called homotopy
      Cohen-Macauley in the literature.
    }
    \begin{example}[The Solomon-Tits Theorem]
      Every spherical building (i.e., a building with
      finite Weil group) is a spherical simplicial complex. Since
      all links in spherical buildings are spherical buildings,
      we see that spherical buildings are Cohen-Macauley.
    \end{example}
    \begin{example}[{\cite[Theorem 12.4]{Quillen:1978}}]
      Let $\ThePrime$ be a prime number and assume that the field
      $\TheField$ has characteristic
      $\neq\ThePrime$ and contains
      a $\ThePrime^{\text{th}}$ root of unity. Then the poset
      $\ElementaryAbelianSubgroupPosetOf[\ThePrime]{\GlOf[\TheRank]{\TheField}}$
      of non-trivial
      elementary Abelian subgroups of $\GlOf[\TheRank]{\TheField}$
       is Cohen-Macauley of dimension
      $\TheRank-\One$.
    \end{example}
    
    We call a $\Delta$-complex Cohen-Macauley if its
    associated poset of simplices is Cohen-Macauley. Equivalently,
    a $\Delta$-complex is Cohen-Macauley if its barycentric subdivision
    is Cohen-Macauley. Note that the barycentric subdivision of a
    Cohen-Macauley simplicial complex is Cohen-Macauley, too. Thus
    for $\Delta$-complexes that are already simplicial, the two notions
    of being Cohen-Macauley coincide.
    \begin{observation}\label{delta:cohen-macauley}
      We call a $\Delta$-complex \notion{strict} if every closed simplex
      has an injective attaching map, i.e., no faces of an individual simplex
      are identified. In this case, the associated poset locally looks like
      the poset of a simplicial complex: the boundary of any element is isomorphic
      to the poset of strict subsets of a finite set. It follows that
      a strict $\Delta$-complex is Cohen-Macauley if it is spherical and
      has spherical links only.\qed
    \end{observation}
    All complexes discussed below are strict $\Delta$-complexes.

  \section{The Chessboard Complex and its Braided Version}
    The \notion{$\TheBotNumber\times\TheTopNumber$-chessboard complex}
    \(
      \ChessboardComplex[\TheBotNumber][\TheTopNumber]
    \)
    is
    the simplicial complex whose vertex set is the set of squares of an
    $\TheBotNumber\times\TheTopNumber$ chessboard and whose simplices are
    configurations of non-threatening rooks on said chessboard. Thus, formally,
    this simplicial complex is the collection of those subsets
    \(
      \TheSimplex
      \subseteq
      \SetOf{\One,\ldots,\TheBotNumber}
      \crossprod
      \SetOf{\One,\ldots,\TheTopNumber}
    \)
    such that the projections
    \(
      \TheProjection[\TheIndex]
    \)
    ($\TheIndex=\One,\Two$) restrict to injective maps on $\TheSimplex$: we have at
    most one rook in each row and each column of the chessboard.
    
    The complexes $\ChessboardComplex[\TheBotNumber][\TheTopNumber]$ have
    been studied intensively. In particular, a good deal is known about their
    connectivity properties:
    \begin{cit}[{\cite{Bjoerner.Lovasz.Vrecica.Zivaljevic:1994}}]\label{braid:chessboard}
      Put
      \(
        \TheDummyNumber
        :=
          \MinOf{
            \TheBotNumber,
            \TheTopNumber,
            \lfloor
              \frac{\TheBotNumber+\TheTopNumber+\One}{\Three}
            \rfloor
          }
        .
      \)
      Then the chessboard complex
      $\ChessboardComplex[\TheBotNumber][\TheTopNumber]$ is
      $\Parentheses{\TheDummyNumber-\Two}$-connected.
      In fact, the $\Parentheses{\TheDummyNumber-\One}$-skeleton
      of $\ChessboardComplex[\TheBotNumber][\TheTopNumber]$ is Cohen-Macauley.

      In particular, if
      \(
          \MinOf{
            \TheBotNumber,
            \TheTopNumber
          }
          \leq
          \lfloor
            \frac{\TheBotNumber+\TheTopNumber+\One}{\Three}
          \rfloor
      \)
      then the complex $\ChessboardComplex[\TheBotNumber][\TheTopNumber]$ is
      Cohen-Macauley.
    \end{cit}

    There is another description of
    $\ChessboardComplex[\TheBotNumber][\TheTopNumber]$
    in terms of matchings in the complete
    $\Pair{\TheBotNumber}{\TheTopNumber}$-bipartite graph. Recall
    that a \notion{matching} in a graph is a subgraph that consists of
    disjoint edges, i.e., every vertex is contained in at most one edge
    of the subgraph. Table~\ref{braid:rooks_and_matchings} illustrates how a non-threatening
    configuration of rooks and a matching in the complete
    bipartite graph represent the same subset of
    \(
      \SetOf{\One,\ldots,\TheBotNumber}
      \crossprod
      \SetOf{\One,\ldots,\TheTopNumber}
      .
    \)
    \begin{table}
      \begin{center}
        \begin{tabular}{@{}c@{$\,\,\longleftrightarrow\,\,$}c@{}}
          \begingroup\chessfont
          \def\CSpace{\kern2pt}
          \begin{tabular}{|@{\CSpace}c@{\CSpace}|@{\CSpace}c@{\CSpace}|@{\CSpace}c@{\CSpace}|@{\CSpace}c@{\CSpace}|@{\CSpace}c@{\CSpace}|@{\CSpace}c@{\CSpace}|}
          \hline
            0&0&0&R&0&0\\
          \hline
            0&0&0&0&0&0\\
          \hline
            R&0&0&0&0&0\\
          \hline
            0&0&R&0&0&0\\
          \hline
          \end{tabular}
          \endgroup
          &
          \begin{tabular}{@{}c@{}}
           \fbox{\includegraphics{braids.3}{}{}}
          \end{tabular}
          \\
        \end{tabular}
      \end{center}
      \caption{\label{braid:rooks_and_matchings}Two equivalent ways of representing
        the simplex
        \(
          \SetOf{
            \Pair{1}{3};
            \Pair{2}{1};
            \Pair{4}{4}
          }
        \)}
    \end{table}

    The matching picture suggest the following construction: Embed the vertex
    set of the $\Pair{\TheBotNumber}{\TheTopNumber}$-bipartite graph into the
    boundary of a cube so that the 
    $\TheBotNumber$ blue vertices lie in the bottom square and the
    $\TheTopNumber$ red vertices lie in the top square. A \notion{partial braid}
    is a braid running vertically through the cube all of whose strands connect
    a red to a blue vertex (see table~\ref{table:partial_braid}).
    Of course, two braids are equal if one can be deformed
    into the other by an ambient homotopy fixing the boundary of the cube pointwise.
    \begin{table}
      \begin{center}
      \begin{tabular}{cc}
        \begin{tabular}{c}
          \includegraphics[scale=1]{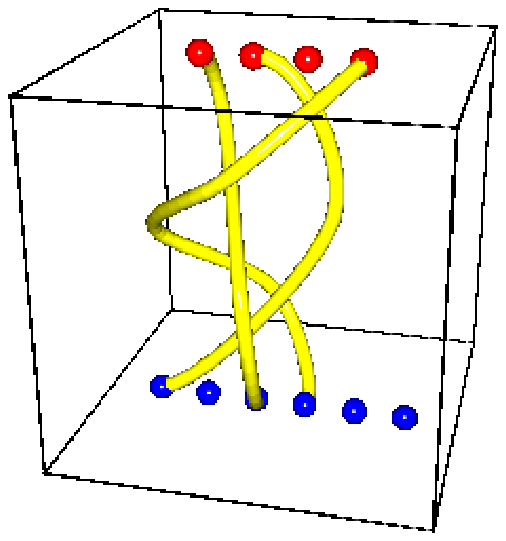}{}{}
        \end{tabular}
        &
        \IncludeMathDiagram{braids.2}
      \end{tabular}
      \end{center}
      \caption{\label{table:partial_braid}A partial braid and its diagram}
    \end{table}
    
    The set $\TheBraidsPoset[\TheBotNumber][\TheTopNumber]$
    of partial braids carries a natural poset structure: The face relation
    $\isstrictfaceof$ is defined by deleting strands. We call the poset thus
    defined the \notion{braided chessboard poset}. This set also is a
    strict $\Delta$-complex in an obvious way: a partial braid on $\TheDim+\One$
    strands is a $\TheDim$-simplex. The poset
    $\TheBraidsPoset[\TheBotNumber][\TheTopNumber]$ is just the
    poset of cells in this complex, and we will silently identify the
    braided chessboard poset
    and the braided chessboard complex.
    Note that the braided chessboard complex
    $\TheBraidsPoset[\TheBotNumber][\TheTopNumber]$ is
    not a simplicial complex:
    all cells have the shape of a simplex, but such a simplex is not determined
    uniquely by its set of vertices.

    We have a natural projection
    \[
      \TheProjection \mapcolon
      \TheBraidsPoset[\TheBotNumber][\TheTopNumber]
      \rightarrow
      \ChessboardPoset[\TheBotNumber][\TheTopNumber]
    \]
    from the braided chessboard poset
    to the poset $\ChessboardPoset[\TheBotNumber][\TheTopNumber]$
    of simplices in $\ChessboardComplex[\TheBotNumber][\TheTopNumber]$ which is
    given by viewing the braid as a matching. This projection is a height-preserving
    morphism of posets. Table~\ref{table:projection} illustrates that the
    face relation and the projection are compatible.
    \begin{table}
      \[
      \begin{array}{ccc}
        \IncludeMathDiagram{braids.4} & \isstrictfaceof & \IncludeMathDiagram{braids.2} \\
        \downarrow                    &           & \downarrow \\
        \IncludeMathDiagram{braids.5} & \subset   & \IncludeMathDiagram{braids.3} \\
      \end{array}
      \]
      \caption{\label{table:projection}Face relation and projection}
    \end{table}
    
    We remark that the braid groups $\BraidGroup[\TheBotNumber]$
    and $\BraidGroup[\TheTopNumber]$ act from opposite sides on the
    braided chessboard complex $\TheBraidsPoset[\TheBotNumber][\TheTopNumber]$ just as
    the symmetric groups $\SymmetricGroup[\TheBotNumber]$ and
    $\SymmetricGroup[\TheTopNumber]$ act on the chessboard complex.

    \begin{observation}
      A very useful property of chessboard complexes is that links in
      chessboard complexes are themselves chessboard complexes of smaller size:
      the link of a simplex of dimension $\TheDim$ in
      \(
        \ChessboardComplex[\TheBotNumber][\TheTopNumber]
      \)
      is isomorphic to the chessboard complex
      \(
        \ChessboardComplex[\TheBotNumber-\TheDim-\One][\TheTopNumber-\TheDim-\One]
        .
      \)
      \qed
    \end{observation}
    Unfortunately, this does not hold for the restricted class of braided
    chessboard complexes that we have defined so far. To make our arguments
    amenable to induction, we therefore generalize our construction slightly:
    Let $\TheCutOutGraph$ be a graph embedded into the cube avoiding the red
    and blue vertices. A \notion{partial braid relative to $\TheCutOutGraph$}
    is a braid running vertically through the cube not meeting a small, closed
    regular neighborhood of $\TheCutOutGraph$ (see table~\ref{table:relative_braid}).
    We consider relative partial braids
    equal if we can deform one into the other by an ambient isotopy that leaves
    the boundary of the cube and the regular neighborhood of $\TheCutOutGraph$
    fixed pointwise.
    \begin{table}
      \begin{center}
        \includegraphics[scale=1]{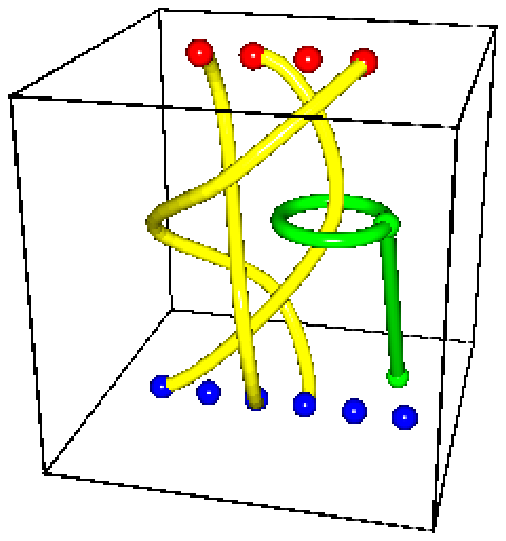}{}{}
      \end{center}
      \caption{\label{table:relative_braid}A partial braid relative to a graph}
    \end{table}
    By
    \(
      \TheBraidsPosetOf[\TheBotNumber][\TheTopNumber]{\TheCutOutGraph}
      ,
    \)
    we denote the poset ($\Delta$-complex)
    of partial braids relative to $\TheCutOutGraph$.
    Note that we still have a natural projection
    \[
      \TheProjection \mapcolon
      \TheBraidsPosetOf[\TheBotNumber][\TheTopNumber]{\TheCutOutGraph}
      \rightarrow
      \ChessboardPoset[\TheBotNumber][\TheTopNumber]
      .
    \]
    Now, we can describe links rather easily:
    \begin{observation}
      Let $\TheSimplex$ be a simplex of dimension $\TheDim$ in
      \(
        \TheBraidsPosetOf[\TheBotNumber][\TheTopNumber]{\TheCutOutGraph}
        .
      \)
      Then $\TheSimplex$ is itself realized as a graph embedded in the cube.
      Its link is isomorphic to
      \(
        \TheBraidsPosetOf[\TheBotNumber-\TheDim-\One][\TheTopNumber-\TheDim-\One]{
          \TheSimplex\union\TheCutOutGraph
        }
        .
      \)
      \qed
    \end{observation}

    Our goal is to understand the connectivity of the posets
    \(
      \TheBraidsPosetOf[\TheBotNumber][\TheTopNumber]{\TheCutOutGraph}
    \)  
    and, in order to derive their connectivity properties,
    we will first study the fibers of the projection
    \(
      \TheProjection \mapcolon
      \TheBraidsPosetOf[\TheBotNumber][\TheTopNumber]{\TheCutOutGraph}
      \rightarrow
      \ChessboardPoset[\TheBotNumber][\TheTopNumber]
      .
    \)

    \begin{prop}\label{braid:fiber_connectivity}
      Let $\AltSimplex$ be a simplex of dimension $\TheDim$ in the
      chessboard complex $\ChessboardComplex[\TheBotNumber][\TheTopNumber]$
      Then the \notion{closed fiber} over $\AltSimplex$
      \[
        \FiberOver{\AltSimplex}
        :=
        \TheProjectionOf[][-\One]{\ClosureOf{\AltSimplex}}
      \]
      is Cohen-Macauley of dimension $\Parentheses{\TheDim-\One}$.
    \end{prop}
    The proof is by induction. The case $\TheDim=\Zero$ just states
    that the fiber over a vertex is non-empty and discrete. Both claims
    are clearly true: any graph consisting of a single edge can be lifted,
    and any two such lifts
    are not joined by higher dimensional material.
    In the proofs
    of the following two lemmas, we will assume that the proposition holds
    for all simplices of dimension strictly less than $\TheDim$ and all
    graphs $\TheCutOutGraph$. In particular,
    we will assume the proposition for all strict faces of the simplex $\AltSimplex$.
    \begin{lemma}
      Let $\TheSimplex$ be a strict face of the simplex $\AltSimplex$
      in the chessboard complex $\ChessboardComplex[\TheBotNumber][\TheTopNumber]$.
      Then the inclusion
      \(
        \FiberOver{\TheSimplex}
        \monorightarrow
        \FiberOver{\AltSimplex}
      \)
      is trivial in homotopy.
    \end{lemma}
    \begin{proof}
      Fix a sphere $\TheSphere$ inside $\FiberOver{\TheSimplex}$. We have to
      show that this sphere can be contracted inside $\FiberOver{\AltSimplex}$.
      
      First note that, by compactness, the sphere $\TheSphere$ involves only
      finitely many cells in $\FiberOver{\TheSimplex}$. Each of these cells is
      represented by some partial braid in
      \(
        \TheBraidsPosetOf[\TheBotNumber][\TheTopNumber]{\TheCutOutGraph}
      \)
      that lies above $\ClosureOf{\TheSimplex}$. 
      Let $\TheVertex$ be a vertex in $\AltSimplex\setminus\TheSimplex$.
      This vertex represents an edge in the
      matching $\AltSimplex$ connecting a red and a blue vertex.
      In the cube, we can
      find a strand $\TheStrand$
      connecting this red vertex to this blue vertex, a strand
      that does not braid
      with any of the partial braids used in the sphere $\TheSphere$:
      there is only a finite number of them and we can always work around finitely
      many partial braids by staying sufficiently close to the boundary of the
      cube.
      The strand $\TheStrand$, therefore, can be added to all the partial braids
      in a way that is compatible with face relations. Thus,
      this strand $\TheStrand$ represents a vertex
      in $\FiberOver{\AltSimplex}$ that serves
      as a cone point for a contraction of the sphere $\TheSphere$.
    \end{proof}
    \begin{rem}
      It is not really important that the strand $\TheStrand$ does not braid with
      any of the finitely many given partial braids used to build $\TheSphere$.
      This is just the easiest way to ensure that $\TheStrand$ braids with all
      the given braids in a consistent way: Since braids allow for deformation,
      a common face to two partial braids can look very different in its cofaces.
      We have to add $\TheStrand$ in such a way that face relations are preserved.
    \end{rem}
    We now turn to links inside the closed fiber $\FiberOver{\AltSimplex}$.
    The slogan is: links in a fiber are fibers in a link:
    \begin{observation}\label{braid:fiber_link}
      Let $\TheVertex$ be a vertex in $\AltSimplex$.
      Note that the fiber over $\TheVertex$ is a discrete set: every element
      is one way of embedding a single strand in the complement of $\TheCutOutGraph$.
      Fix a strand $\TheStrand$ representing a vertex in the fiber
      $\FiberOver{\TheVertex}\subset\FiberOver{\AltSimplex}$.
      The projection
      \(
        \TheProjection
        \mapcolon
        \TheBraidsPosetOf[\TheBotNumber][\TheTopNumber]{\TheCutOutGraph}
        \rightarrow
        \ChessboardPoset[\TheBotNumber][\TheTopNumber]
      \)
      is a strictly monotonic poset-map (i.e., it does not collapse cells).
      It restricts to the link of $\TheStrand$ as follows:
      \[
        \begin{array}{ccc}
          \LinkOf[{
            \TheBraidsPosetOf[\TheBotNumber][\TheTopNumber]{\TheCutOutGraph}
          }]{
            \TheStrand
          }
          &=&
          \TheBraidsPosetOf[\TheBotNumber-\One][\TheTopNumber-\One]{
            \TheCutOutGraph
            \union
            \TheStrand
          }
          \\
          \Big\downarrow\vcenter{
            \rlap{$\TheProjection[\TheCutOutGraph]$}
          }
          & &
          \Big\downarrow\vcenter{
            \rlap{$\TheProjection[\TheCutOutGraph\union\TheStrand]$}
          }
          \\
          \LinkOf[{\ChessboardPoset[\TheBotNumber][\TheTopNumber]}]{\TheVertex}
          &=&
          \ChessboardPoset[\TheBotNumber-\One][\TheTopNumber-\One]
        \end{array}
      \]
      Moreover, the link of $\TheStrand$ inside $\FiberOver{\AltSimplex}$ is
      the part in the link that maps to $\AltSimplex\setminus\SetOf{\TheVertex}$.
      Thus:
      \[
        \LinkOf[\FiberOver{\AltSimplex}]{\TheStrand}
        =
        \FiberOver[{
          \TheProjection[\TheCutOutGraph\union\TheStrand]
        }]{\AltSimplex\setminus\SetOf{\TheVertex}}
        .
        \qed
      \]
    \end{observation}
    \begin{lemma}
      Let $\TheSimplex$ be a codimension $\One$ face of the simplex $\AltSimplex$
      in the chessboard complex $\ChessboardComplex[\TheBotNumber][\TheTopNumber]$.
      Then the inclusion
      \(
        \FiberOver{\TheSimplex}
        \monorightarrow
        \FiberOver{\AltSimplex}
      \)
      induces epimorphisms
      \(
        \HomotopyGroupOf[\TheIndex]{
          \FiberOver{\TheSimplex}
        }
        \epirightarrow
        \HomotopyGroupOf[\TheIndex]{
          \FiberOver{\AltSimplex}
        }
      \)
      in homotopy
      for $\TheIndex < \DimOf{\AltSimplex}$.
    \end{lemma}
    \begin{proof}
      Let $\TheVertex$ be the vertex opposite to $\TheSimplex$ in $\AltSimplex$
      and let $\TheStrand$ be a strand representing a vertex
      in $\FiberOver{\TheVertex}$. We have
      just observed that the relative link
      $\LinkOf[\FiberOver{\AltSimplex}]{\TheStrand}$
      is the fiber above $\TheSimplex$ in the complex
      \(
        \TheBraidsPosetOf[\TheBotNumber-\One][\TheTopNumber-\One]{
          \TheCutOutGraph
          \union
          \TheStrand          
        }
        .
      \)
      Since $\DimOf{\TheSimplex}<\DimOf{\AltSimplex}$ we can apply
      Proposition~\ref{braid:fiber_connectivity} by induction. It follows
      that the relative link of $\TheStrand$ is
      $\Parentheses{\DimOf{\AltSimplex}-\Two}$-connected.
      Thus, for $\TheIndex\leq\DimOf{\AltSimplex}-\One$,
      every $\TheIndex$-sphere that passes
      through $\TheStrand$ can be homotoped off the vertex $\TheStrand$.
      We can do this simultaneously for all vertices above $\TheVertex$ and push
      any $\TheIndex$-sphere into the fiber $\FiberOver{\TheSimplex}$. Therefore,
      this fiber carries all of $\HomotopyGroupOf[\TheIndex]{\FiberOver{\AltSimplex}}$.
    \end{proof}
    \begin{rem}
      A more formal proof can be based on combinatorial Morse theory for piecewise
      Euclidean complexes: Let
      \(
        \TheFunction \mapcolon \AltSimplex \rightarrow
        \ClosedInterval{\Zero}{\One}
      \)
      be the affine map sending $\TheVertex$ to $\One$ and $\TheSimplex$ to
      $\Zero$. The composition $\TheFunction\compose\TheProjection$ is a
      Morse function on $\FiberOver{\AltSimplex}$ in the sense of
      \cite{Bestvina.Brady:1997}. The preceding argument
      establishes that the descending links are
      $\Parentheses{\DimOf{\AltSimplex}-\Two}$-connected. Now the lemma follows
      from \cite[Lemma~2.5 and Corollary~2.6]{Bestvina.Brady:1997}.
    \end{rem}
    \begin{proof}[of Proposition~\ref{braid:fiber_connectivity} (finish)]
      The preceding two lemmas state that the map
      \[
        \HomotopyGroupOf[\TheIndex]{
          \FiberOver{\TheSimplex}
        }
        \epirightarrow
        \HomotopyGroupOf[\TheIndex]{
          \FiberOver{\AltSimplex}
        }
      \]
      is trivial and onto for $\TheIndex<\DimOf{\AltSimplex}$. Thus we
      know that the fiber $\FiberOver{\AltSimplex}$ is
      $\Parentheses{\TheDim-\One}$-connected. Thus fibers are spherical.

      Invoking Observation~\ref{braid:fiber_link} again, we conclude that
      all links of cells in $\FiberOver{\AltSimplex}$ are spherical, too.
      Since $\FiberOver{\AltSimplex}$ is a strict $\Delta$-complex, we
      infer that it is Cohen-Macauley by Observation~\ref{delta:cohen-macauley}.
    \end{proof}

    Now that we understand the fibers of the projection, we can apply
    the tools provided by D.~Quillen:
    \begin{cit}[{\cite[Corollary~9.7]{Quillen:1978}}]\label{braid:quillen_cm}
      Suppose
      \(
        \ThePosetMap
        \mapcolon
        \ThePoset \rightarrow \AltPoset
      \)
      is a strictly increasing
      morphism of posets. Assume that $\AltPoset$ is
      Cohen-Macauley of dimension $\TheDim$ and
      that for every $\AltPosetElement\in\AltPoset$, the preimage
      \(
        \PreImage{\ThePosetMap}{\ClosureOf{\AltPosetElement}}
        =
        \SetOf[
          \ThePosetElement\in\ThePoset
        ]{
          \ThePosetMapOf{\ThePosetElement}\isfaceof\AltPosetElement
        }
      \)
      of the closure
      $\ClosureOf{\AltPosetElement}$ is Cohen-Macauley of
      dimension $\HeightOf{\AltPosetElement}$. Then $\ThePoset$ is Cohen-Macauley of
      dimension $\TheDim$.
    \end{cit}
    \begin{theorem}
      Put
      \(
        \TheDummyNumber
        :=
          \MinOf{
            \TheBotNumber,
            \TheTopNumber,
            \lfloor
              \frac{\TheBotNumber+\TheTopNumber+\One}{\Three}
            \rfloor
          }
        .
      \)
      Then, for any graph $\TheCutOutGraph$ in the cube, the braided chessboard
      complex
      \(
        \TheBraidsPosetOf[\TheBotNumber][\TheTopNumber]{
          \TheCutOutGraph
        }
      \)
      has a Cohen-Macauley $\Parentheses{\TheDummyNumber-\One}$-skeleton. In
      particular, the complex is $\Parentheses{\TheDummyNumber-\Two}$-connected;
      and if
      \(
          \MinOf{
            \TheBotNumber,
            \TheTopNumber
          }
          \leq
          \lfloor
            \frac{\TheBotNumber+\TheTopNumber+\One}{\Three}
          \rfloor
      \)
      then $\TheBraidsPosetOf[\TheBotNumber][\TheTopNumber]{\TheCutOutGraph}$ is
      Cohen-Macauley.
    \end{theorem}
    \begin{proof}
      Since the projection
      \(
        \TheProjection \mapcolon
        \TheBraidsPosetOf[\TheBotNumber][\TheTopNumber]{\TheCutOutGraph}
        \rightarrow
        \ChessboardPoset[\TheBotNumber][\TheTopNumber]
      \)
      does not crush cells, the $\Parentheses{\TheDummyNumber-\One}$-skeleton
      of the braided chessboard complex is the preimage of the
      $\Parentheses{\TheDummyNumber-\One}$-skeleton of the chessboard
      complex $\ChessboardComplex[\TheBotNumber][\TheTopNumber]$, which
      is Cohen-Macauley by Citation~\ref{braid:chessboard}.
      
      By Proposition~\ref{braid:fiber_connectivity}, the projection
      \(
        \TheProjection
      \)
      satisfies the hypotheses of Citation~\ref{braid:quillen_cm}. The theorem follows
      immediately.
    \end{proof}

  \section{Complexes Based on Collections of Edges}
    Let $\TheVertexSet$ be a finite set, and let $\TheGraphFam$ be a family of
    graphs sharing $\TheVertexSet$ as their vertex sets.
    Suppose that $\TheGraphFam$
    is \notion{subgraph-closed}, i.e., if $\TheGraph\in\TheGraphFam$
    then every subgraph of $\TheGraph$ also belongs to $\TheGraphFam$.
    The \notion{graph poset}
    \begin{notationlist}
      \item
        \(
          \GraphPosetOf{\TheGraphFam}
        \)
        induced by $\TheGraphFam$ is the poset
        of non-empty graphs in $\TheGraphFam$
        ordered by inclusion. The \notion{graph complex}
      \item
        \(
          \GraphComplexOf{\TheGraphFam}
        \)
        induced by $\TheGraphFam$ is
        the simplicial complex whose $\TheDim$-simplices are those graphs in
        $\TheGraphFam$ containing precisely $\TheDim+\One$ edges. The graph
        poset is the poset of cells in this complex.
    \end{notationlist}
    \begin{example}[The Complex of Not \boldmath$\TheIndex$-Connected Graphs]
      Consider the family of non-$\TheIndex$-connected
      simplicial graphs over the vertex set $\TheVertexSet$.
      The corresponding graph complexes have been studied in
      \cite{Babson.Bjoerner.Linusson.Shareshian.Welker:1999}. In particular,
      it is shown that the complex of not $\Two$-connected graphs on
      $\TheVertexSet$ is homotopy equivalent to a wedge of
      $\FactorialOf{\Parentheses{\CardOf{\TheVertexSet}-\Two}}$
      spheres of dimension $\Two\CardOf{\TheVertexSet}-\Five$.
    \end{example}
    \begin{example}[The Forest Complex]
      Let $\TheGraph$ be a fixed graph over the
      vertex set $\TheVertexSet$ with $\TheNumberOfComponents$
      components, and let
      $\TheGraphFam$ be the family of forests
      in $\TheGraph$, i.e., circle-free subgraphs
      of $\TheGraph$.
      Then the complex
      \[
        \ForestComplexOf{\TheGraph} :=
        \GraphComplexOf{\TheGraphFam}
      \]
      is the \notion{complex of forest} in $\TheGraph$.
    \end{example}
    \begin{prop}
      The forest complex $\ForestComplexOf{\TheGraph}$ is Cohen-Macauley of dimension
      $\CardOf{\TheVertexSet}-\TheNumberOfComponents-\One$.
    \end{prop}
      This was proved independently
      by several people. The earliest source, I am aware of is the thesis of
      J.S.~Provan \cite{Provan:1977}. The forest complex is the independence complex
      of a matroid and hence shellable by \cite[Theorem 7.3.3]{Bjoerner:1992}.
      For those who are scared by matroids and shellability, we include a down
      to earth proof base on the version given in \cite[Proposition~2.2]{Vogtmann:1990}.
      \begin{proof}
        Every spanning forest of $\TheGraph$ contains precisely
        $\CardOf{\TheVertexSet}-\TheNumberOfComponents$ edges. Thus each
        maximal simplex in 
        $\ForestComplexOf{\TheGraph}$ has
        dimension $\CardOf{\TheVertexSet}-\TheNumberOfComponents-\One$.

        A simplex $\TheSimplex$ in $\ForestComplexOf{\TheGraph}$ is a sub-forest
        of $\TheGraph$. Collapsing this sub-forest yields a new graph
        $\TheGraph\rmod\TheSimplex$ that has the same number of connected
        components. However, each edge in the forest $\TheSimplex$ connects two
        vertices, whence crushing this edge reduces the number of vertices by
        $\One$. It follows that $\ForestComplexOf{\TheGraph\rmod\TheSimplex}$
        has dimension
        $\DimOf{\ForestComplexOf{\TheGraph}}-\DimOf{\TheSimplex}-\One$. The link
        of $\TheSimplex$ in $\ForestComplexOf{\TheGraph}$ is
        isomorphic to $\ForestComplexOf{\TheGraph\rmod\TheSimplex}$. Thus
        the complex will be Cohen-Macauley, provided that
        $\ForestComplexOf{\TheGraph}$ is spherical for all graphs $\TheGraph$.
        Since we already established the dimension of $\ForestComplexOf{\TheGraph}$,
        it remains to show that $\ForestComplexOf{\TheGraph}$ is
        $\Parentheses{\CardOf{\TheVertexSet}-\TheNumberOfComponents-\Two}$-connected.

        Let $\TheEdge$ represent a vertex in $\ForestComplexOf{\TheGraph}$,
        i.e, $\TheEdge$ is a non-loop edge in $\TheGraph$. If $\TheEdge$ is
        a separating edge, then it serves as a cone point in
        $\ForestComplexOf{\TheGraph}$, in which case the forest complex
        is contractible and a fortiori
        $\Parentheses{\CardOf{\TheVertexSet}-\TheNumberOfComponents-\Two}$-connected.

        If $\TheEdge$ is non-separating, we can remove $\TheEdge$ without increasing
        the number of components. Thus the graph $\TheGraph\setminus\TheEdge$ has
        $\CardOf{\TheVertexSet}$ vertices and $\TheNumberOfComponents$ components.
        By induction on the number of edges, we may assume that
        $\ForestComplexOf{\TheGraph\setminus\TheEdge}$ is spherical of dimension
        $\Parentheses{\CardOf{\TheVertexSet}-\TheNumberOfComponents-\Two}$. From
        this subcomplex, we obtain $\ForestComplexOf{\TheGraph}$ by coning of
        the link of $\TheEdge$, which is isomorphic to
        $\ForestComplexOf{\TheGraph\rmod\TheEdge}$, which is spherical of
        dimension
        $\Parentheses{\CardOf{\TheVertexSet}-\One-\TheNumberOfComponents-\Two}$
        again by induction on the number of edges. It follows that
        $\ForestComplexOf{\TheGraph}$ is spherical of dimension
        $\Parentheses{\CardOf{\TheVertexSet}-\TheNumberOfComponents-\Two}$.
      \end{proof}

    Another family of examples arises as follows: Fix a
    graph $\TheGraph$ with vertex set $\TheVertexSet$,
    and let $\MatchingsOf{\TheGraph}$ be the family of
    subgraphs satisfying the condition that each vertex
    is contained in at most one edge.  Such subgraphs are called
    \notion{matchings} in $\TheGraph$. We denote the
    graph poset associated to the family of matchings
    by $\MatchingPosetOf{\TheGraph}$
    and its graph complex by $\MatchingComplexOf{\TheGraph}$.
    \begin{example}[The Chessboard Complex]
      If $\TheCompleteGraph[\TheBotNumber,\TheTopNumber]$ is the complete bipartite
      graph on $\TheTopNumber$ red and $\TheBotNumber$
      blue vertices, non-empty, edge-disjoint subgraphs correspond
      to partial matchings between the set of
      red vertices and the set of blue vertices. Thus, we
      recover the chessboard complex:
      \[
        \MatchingComplexOf{\TheCompleteGraph[\TheBotNumber,\TheTopNumber]}
        =
        \ChessboardComplex[\TheBotNumber][\TheTopNumber]
        .
      \]
    \end{example}
    We mention that
    \(
      \ChessboardComplex[\TheBotNumber][\TheTopNumber]
    \)
    is sometimes called the matching complex. We prefer,
    however to use this name for the following:
    \begin{example}[The Matching Complex]
      Let $\TheCompleteGraph[\TheNumber]$ be the complete
      graph on $\TheNumber$ vertices. The elements of
      $\MatchingPosetOf{\TheCompleteGraph[\TheNumber]}$ are collections of disjoint
      edges.
      The corresponding graph complex
      $\MatchingComplexOf{\TheCompleteGraph[\TheNumber]}$ is called
      the \notion{matching complex}.

      Some connectivity properties of these complexes are known: Put
      \[
        \TheDummyNumber
        :=
        \left\lfloor
          \frac{\TheNumber+\One}{\Three}
        \right\rfloor
        .
      \]
      Then the
      $\Parentheses{\TheDummyNumber-\One}$-skeleton of
      $\MatchingComplexOf{\TheCompleteGraph[\TheNumber]}$
      is Cohen-Macauley
      \cite[Corollary~4.2]{Bjoerner.Lovasz.Vrecica.Zivaljevic:1994}.
    \end{example}
    
    Finally, we will have a use for the most basic subgraph-closed family:
    \begin{example}[The Simplex]\label{simplex}
      Let $\SubGraphsOf{\TheGraph}$ be the family of
      subgraphs of a given graph $\TheGraph$, then
      $\GraphComplexOf{\SubGraphsOf{\TheGraph}}$ is nothing
      but a big simplex whose vertices are the edges
      in $\TheGraph$. A single simplex is Cohen-Macauley.
    \end{example}

  \section{The Tangling Construction}
    In the case of the chessboard complex, the underlying graph was bipartite.
    Thus, we could put the two kinds of vertices into opposite faces, top and
    bottom, of the cube
    and require that strands pass through the cube vertically. In general, we
    cannot arrange for this. Thus we will replace braids by \notion{tangles}
    to make the construction applicable to the complexes discussed above.

    \begin{table}
      \begin{center}
        \begin{tabular}{c}
          \(
            \begin{array}{c}\includegraphics{braids.10}{}{}\end{array}
            \qquad
            \begin{array}{c}\includegraphics[scale=0.75]{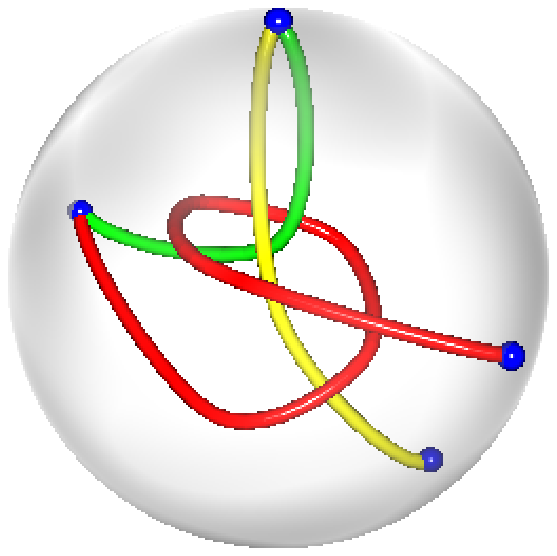}{}{}\end{array}
            \qquad
            \begin{array}{c}\includegraphics[scale=0.75]{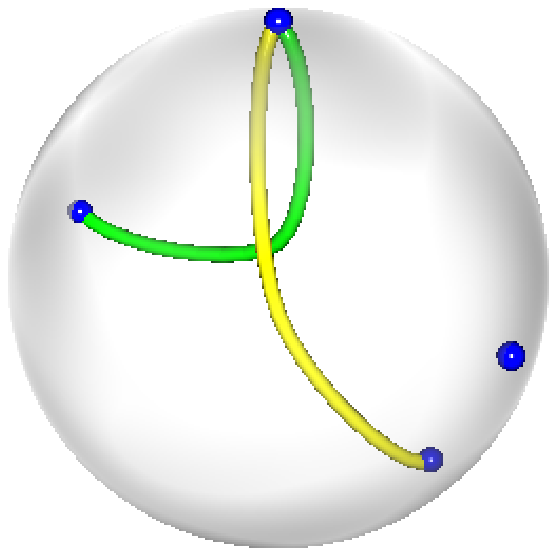}{}{}\end{array}
          \)
        \end{tabular}
      \end{center}
      \caption{\label{tangle:lift}A graph, a lift, and a face}
    \end{table}
    Let $\TheGraphFam$ be a family of graphs over the vertex set $\TheVertexSet$
    and assume $\TheGraphFam$ is subgraph-closed.
    Let $\TheBall$ be a closed $\Three$-ball. Chose an
    embedding of $\TheVertexSet$ into its boundary
    sphere $\TheSphere:=\BoundaryOf{\TheBall}\subset\TheBall$. A
    \notion{lift} of a graph $\TheGraph\in\TheGraphFam$
    is an embedding of $\TheGraph$ into $\TheBall$ that
    (a) extends the embedding of $\TheVertexSet$ and
    (b) maps interior points of edges to
    interior of $\TheBall$ such that
    (c) each edges lifts to an un-knotted
    curve in $\TheBall$ (see table~\ref{tangle:lift}).
    We will also call these lifts \notion{tangles}. The term lift will be
    preferred when we want to stress the relation to the underlying graph in
    $\TheGraphFam$, whereas the term tangle emphasizes the geometric structure
    of the embedded graph upstairs in the $\Three$-ball.
    We consider two tangles equal if there is an ambient
    homotopy between the two fixing the
    boundary sphere of $\TheBall$ (and, therefore, in
    particular the set $\TheVertexSet$) pointwise.
    The set of equivalence classes of tangles forms the
    \notion{tangle poset}
    \begin{notationlist}
      \item
        \(
          \TanglePosetOf{\TheGraphFam}
        \)
        where the face relation is given by
        deletion of strands: a \notion{strand} is the lift of an
        edge. The \notion{tangle complex}
      \item
        \(
          \TangleComplexOf{\TheGraphFam}
        \)
        is the $\Delta$-complex whose $\TheDim$-simplices are indexed
        by tangles with $\TheDim+\One$ strands.
    \end{notationlist}
    
    As we did with the braided chessboard complex, we will generalize this
    construction by allowing that a given embedded graph $\TheCutOutGraph$ be
    removed from the $\Three$-ball from the start. This way, we will make sure
    that the class of $\Delta$-complexes we define is closed with respect
    to taking links: Let $\TheCutOutGraph$ be
    a graph embedded in $\TheBall$. An \notion{$\TheCutOutGraph$-lift}
    of $\TheGraph$ (or a \notion{$\TheCutOutGraph$-tangle}) is an embedding
    of $\TheGraph$ into $\TheBall$ satisfying the conditions (a) to (c) above
    such that the interiors
    of edges do not meet a fixed regular neighborhood of $\TheCutOutGraph$. Again,
    two $\TheCutOutGraph$-tangles are equivalent if there is an
    ambient homotopy from one to the other fixing the boundary sphere and
    the regular neighborhood of $\TheCutOutGraph$ pointwise.
    Equivalence classes of $\TheCutOutGraph$-tangle form the poset
    \begin{notationlist}
      \item
        \(
          \TanglePosetOf[\TheCutOutGraph]{\TheGraphFam}
          ,
        \)
        the face relation being deletion of strands.
    \end{notationlist}
    \begin{observation}\label{braid:links}
      Let $\TheGraph\in\TheGraphFam$ be a graph
      with $\TheDim+\One$ edges, i.e., an element
      of hight $\TheDim$ in the graph
      poset $\GraphPosetOf{\TheGraphFam}$. The link
      \(
        \LinkOf[\GraphPosetOf{\TheGraphFam}]{\TheGraph}
      \)
      consists of those graphs in $\TheGraphFam$
      that contain $\TheGraph$ as a proper subgraph.
      Removing the edges from $\TheGraph$ yields an isomorphism
      \[
        \LinkOf[\GraphPosetOf{\TheGraphFam}]{\TheGraph}
        \cong
        \GraphPosetOf{\TheGraphFam\graphminus\TheGraph}
      \]
      where $\TheGraphFam\graphminus\TheGraph$ is the
      family of those graphs in $\TheGraphFam$ that
      do not share any edges with $\TheGraph$.

      Similarly,
      for any $\TheCutOutGraph$-lift
      $\TheGraphLift$ of $\TheGraph$, we have an isomorphism
      \[
        \LinkOf[{\TanglePosetOf[\TheCutOutGraph]{\TheGraphFam}}]{\TheGraphLift}
        \cong
        \TanglePosetOf[\TheCutOutGraph\union\TheGraphLift]{\TheGraphFam\graphminus\TheGraph}
        .
      \]
      The isomorphism is not given by removing strands in $\TheGraphLift$, but by
      ``freezing'' them.\qed
    \end{observation}
        
    Our goal is to prove:
    \begin{theorem}\label{braid:main}
      If $\GraphPosetOf{\TheGraphFam}$ is Cohen-Macauley,
      then so is $\TanglePosetOf[\TheCutOutGraph]{\TheGraphFam}$.
    \end{theorem}
    This applies in particular to the forest complex, the chess board complex,
    and to some skeleton in the matching complex.
    \begin{rem}
      The condition (c) above requiring strands to be un-knotted
      can be dropped without affecting the theorem. The proof presented
      here applies to the altered construction without change.
    \end{rem}

    Note that there are canonical, strictly monotonic, hight-preserving projections
    \begin{eqnarray*}
      \TanglePosetOf{\TheGraphFam} &\rightarrow& \GraphPosetOf{\TheGraphFam} \\
      \TangleComplexOf{\TheGraphFam} &\rightarrow& \GraphComplexOf{\TheGraphFam}
    \end{eqnarray*}
    defined by ``ignoring the entanglement''.

    We will closely follow the argument given for the braided chessboard complex.
    Thus, we have to understand fibers over closed simplices.
    \begin{lemma}\label{braid:trivial}
      Let $\TheGraph[\Zero]\subset\TheGraph[\One]$ be a strict inclusion of graphs.
      The induced inclusion
      \[
        \TangleComplexOf[\TheCutOutGraph]{\SubGraphsOf{\TheGraph[\Zero]}}
        \monorightarrow
        \TangleComplexOf[\TheCutOutGraph]{\SubGraphsOf{\TheGraph[\One]}}
      \]
      is trivial in homotopy, i.e., any sphere in
      $\TangleComplexOf[\TheCutOutGraph]{\SubGraphsOf{\TheGraph[\Zero]}}$
      can be crushed
      inside
      $\TangleComplexOf[\TheCutOutGraph]{\SubGraphsOf{\TheGraph[\One]}}$.
    \end{lemma}
    \begin{proof}
      Assume first that $\TheCutOutGraph$ is empty.
      W.l.o.g., we can assume that $\TheGraph[\Zero]$
      is obtained from $\TheGraph[\One]$ by removing
      precisely one edge $\TheEdge$ connecting, say,
      the vertices $\TheVertex$ and $\AltVertex$. We choose
      a path $\ThePath$ inside the boundary sphere $\BoundaryOf{\TheBall}$
      connecting $\TheVertex$ and
      $\AltVertex$. Any sphere in
      $\TangleComplexOf{\SubGraphsOf{\TheGraph[\Zero]}}$
      involves only finitely many
      simplices. By compactness of $\TheCutOutGraph$-lifts,
      we can push the path $\ThePath$ slightly into the
      interior of $\TheBall$ without meeting any
      strands used by the sphere. After pushing it into the ball, the strand
      $\ThePath$ represents a vertex in
      $\TangleComplexOf{\SubGraphsOf{\TheGraph[\One]}}$
      that allows us to cone off the sphere.
      Thus all homotopy of
      $\TangleComplexOf{\SubGraphsOf{\TheGraph[\Zero]}}$
      dies
      in $\TangleComplexOf{\SubGraphsOf{\TheGraph[\One]}}$.

      The argument works as well for non-empty $\TheCutOutGraph$ -- we
      just observe that in pushing $\ThePath$ we can also
      avoid the neighborhood  of $\TheCutOutGraph$ since it is compact.
    \end{proof}
    \begin{prop}\label{tangle:simplex_lift}
      Let $\TheGraph$ be a graph with $\TheDim+\One$ edges.
      Then the poset
      \(
        \TangleComplexOf[\TheCutOutGraph]{\SubGraphsOf{\TheGraph}}
      \)
      is spherical of dimension $\TheDim$.
    \end{prop}
    \begin{proof}
      We use induction on $\TheDim$. The case $\TheDim=\Zero$
      is obvious. So let $\TheEdgeLift$ be
      a strand representing vertex in
      $\TangleComplexOf[\TheCutOutGraph]{\SubGraphsOf{\TheGraph}}$,
      and let $\TheEdge$ be the edge in
      $\TheGraph$ corresponding to $\TheEdgeLift$.
      Note that the complex
      \(
        \TangleComplexOf[\TheCutOutGraph]{\SubGraphsOf{\TheGraph\graphminus\TheEdge}}
      \)
      is a subcomplex of
      \(
        \TangleComplexOf[\TheCutOutGraph]{\SubGraphsOf{\TheGraph}}
        .
      \)
      Moreover, observe that
      \(
        \TangleComplexOf[\TheCutOutGraph]{\SubGraphsOf{\TheGraph}}
      \)
      is obtained from the subcomplex
      \(
        \TangleComplexOf[\TheCutOutGraph]{\SubGraphsOf{\TheGraph\graphminus\TheEdge}}
      \)
      by coning off
      \[
        \LinkOf{\TheEdgeLift}
        \cong
        \TangleComplexOf[\TheCutOutGraph\union\TheEdgeLift]{
          \SubGraphsOf{\TheGraph}\graphminus\TheEdge
        }
        =
        \TangleComplexOf[\TheCutOutGraph\union\TheEdgeLift]{
          \SubGraphsOf{\TheGraph\graphminus\TheEdge}
        }
      \]
      along the canonical map
      \[
        \TangleComplexOf[\TheCutOutGraph\union\TheEdgeLift]{
          \SubGraphsOf{\TheGraph\graphminus\TheEdge}
        }
        \rightarrow
        \TangleComplexOf[\TheCutOutGraph]{
          \SubGraphsOf{\TheGraph\graphminus\TheEdge}
        }
      \]
      given by deleting $\TheEdgeLift$. Both, the link
      \(
        \LinkOf{\TheEdgeLift}
      \)
      and the subcomplex
      \(
        \TangleComplexOf[\TheCutOutGraph]{\SubGraphsOf{\TheGraph\graphminus\TheEdge}}
      \)
      are $\Parentheses{\TheDim-\One}$-spherical by induction. Thus,
      \(
        \TangleComplexOf[\TheCutOutGraph]{\SubGraphsOf{\TheGraph}}
      \)
      is obtained from an $\Parentheses{\TheDim-\One}$-spherical complex by
      conning off an $\Parentheses{\TheDim-\One}$-spherical space. This process
      does not alter homotopy groups in dimensions $\leq\TheDim-\Two$, and it can
      only kill but not introduce homotopy in dimension $\TheDim-\One$. Thus
      \(
        \HomotopyOf[\TheIndex]{
          \TangleComplexOf[\TheCutOutGraph]{\SubGraphsOf{\TheGraph}}
        }
        =
        \Zero
      \)
      for all $\TheIndex<\TheDim-\One$, and
      \[
        \HomotopyOf[\TheDim-\One]{
          \TangleComplexOf[\TheCutOutGraph]{
            \SubGraphsOf{\TheGraph\graphminus\TheEdge}
          }
        }
        \epirightarrow
        \HomotopyOf[\TheDim-\One]{
          \TangleComplexOf[\TheCutOutGraph]{
            \SubGraphsOf{\TheGraph}
          }
        }
      \]
      is onto. However, this map is trivial by Lemma~\ref{braid:trivial}.
    \end{proof}
    \begin{cor}\label{braid:cm}
      Let $\TheGraph$ be a graph with $\TheDim+\One$ edges,
      then $\TangleComplexOf[\TheCutOutGraph]{\SubGraphsOf{\TheGraph}}$ is
      Cohen-Macauley of dimension $\TheDim$.
    \end{cor}
    \begin{proof}
      By Observation~\ref{braid:links},
      the link of a simplex $\TheSimplex$ in
      $\TangleComplexOf[\TheCutOutGraph]{\SubGraphsOf{\TheGraph}}$
      is isomorphic to
      \(
        \TangleComplexOf[\TheCutOutGraph\union\TheSimplex]{
          \SubGraphsOf{\TheGraph'}
        }
      \)
      where $\TheGraph'$ is obtained from $\TheGraph$ by deleting all
      edges that have lifts in $\TheSimplex$.
      Thus, the preceding
      proposition applies to those links them as well, and all
      links in $\TangleComplexOf[\TheCutOutGraph]{\SubGraphsOf{\TheGraph}}$
      are spherical. Thus the strict $\Delta$-complex
      $\TangleComplexOf[\TheCutOutGraph]{\SubGraphsOf{\TheGraph}}$
      is Cohen-Macauley.
    \end{proof}

    \begin{proof}[of Theorem~\ref{braid:main}]
      The claim follows at once from Quillen's
      result~\ref{braid:quillen_cm} and Corollary~\ref{braid:cm} because
      the closed fiber of the projection
      \(
        \TheProjection[\TheCutOutGraph]
        \mapcolon
        \TanglePosetOf[\TheCutOutGraph]{\TheGraphFam}
        \rightarrow
        \GraphPosetOf{\TheGraphFam}
      \)
      above the simplex represented by $\TheGraph\in\TheGraphFam$ is
      isomorphic to
      $\TanglePosetOf[\TheCutOutGraph]{\SubGraphsOf{\TheGraph}}$.
    \end{proof}

  \section{My Motivation: The Complex of Pinched Braids}
    Finally, I would like to present another $\Delta$-complex that also
    projects onto the chessboard complex. It is my struggle with this complex that
    motivated the study of the (seemingly more natural) constructions
    discussed above.

    \begin{table}
      \begin{center}
        \includegraphics{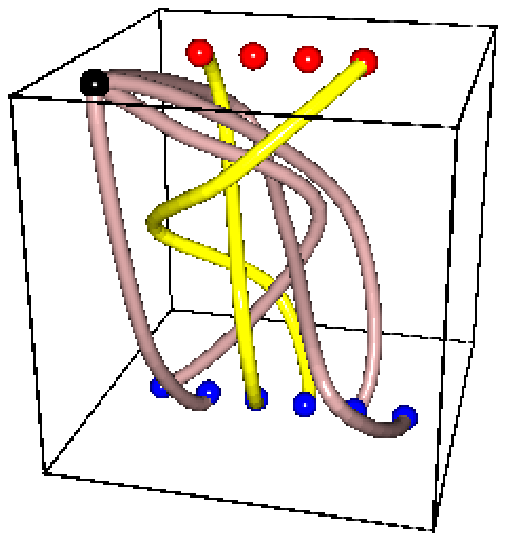}{}{}
      \end{center}
      \caption{\label{pinched:braid}A pinched braid with two regular strands}
    \end{table}
    Fix two numbers $\TheBotNumber$ and $\TheTopNumber$. Embed a row of
    $\TheBotNumber$ blue vertices in the bottom face of the cube labelled
    from left to right by $\TheBotNumber,\TheBotNumber-\One,\ldots,\Two,\One$.
    Embed a row
    of $\TheTopNumber$ red vertices in the top face and
    label them from left to right by $\One,\Two,\ldots,\TheTopNumber-\One,\TheTopNumber$.
    In front of the red row
    add a black vertex in the top face. We think of this vertex as being labelled
    by $\infty$. A \notion{pinched braid} is a collection of disjoint strands
    running vertically through the interior of the
    cube connecting top vertices to bottom vertices
    such that the following conditions are met:
    \begin{enumerate}
      \item
        Every bottom vertex is hit by precisely one strand.
      \item
        Every red vertex is hit by at most one strand. 
      \item
         At least one red vertex is hit by a strand.
    \end{enumerate}
    Strands
    hitting red vertices are called \notion{regular}. The strands
    issuing from $\infty$ are called \notion{singular}.
    The conditions imply, that generically, there will be several singular
    strands, i.e., the black vertex will issue
    multiple strands. We consider two pinched braids equal if they can be deformed
    into one another by ambient homotopies fixing the boundary of the cube
    pointwise. Table~\ref{pinched:braid} shows a pinched braid with two
    regular strands.
    
    \begin{table}
      \[
        \begin{array}{c}
          \includegraphics[scale=0.6]{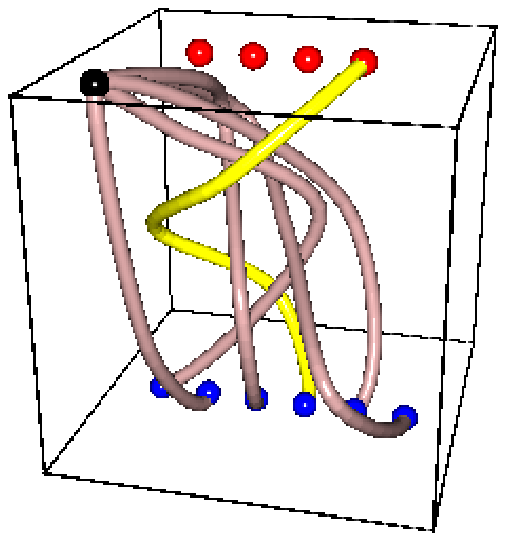}{}{}
        \end{array}
        \isstrictfaceof
        \begin{array}{c}
          \includegraphics[scale=0.6]{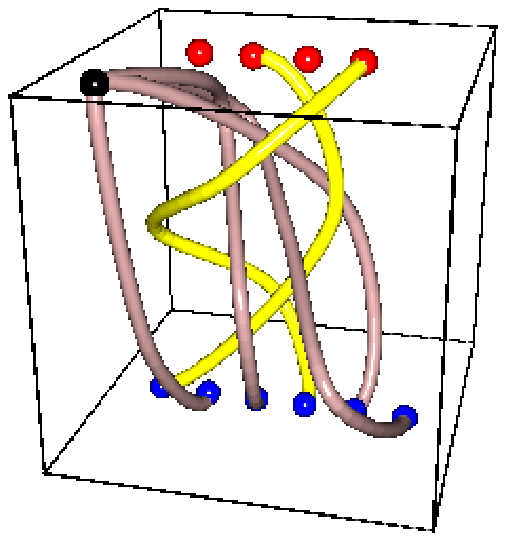}{}{}
        \end{array}
        \isstrictfaceof
        \begin{array}{c}
          \includegraphics[scale=0.6]{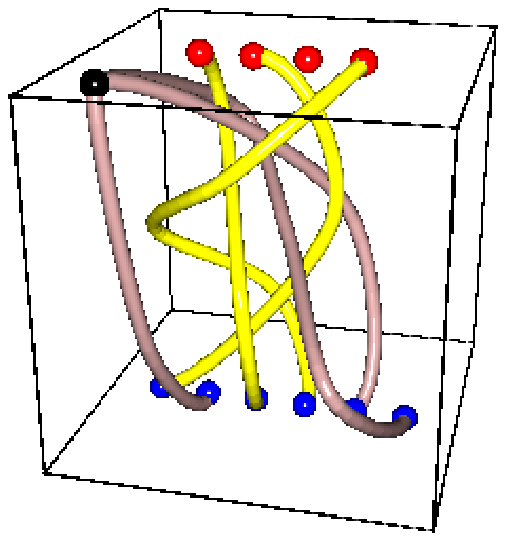}{}{}
        \end{array}
      \]
      \caption{\label{pinched:face_relation}A face chain of length two}
    \end{table}
    The set $\PinchedBraidPoset[\TheBotNumber][\TheTopNumber]$ of all
    pinched braids forms a poset with the following face relation: every
    pinched braid has one immediate face for each regular strand
    obtained by sliding
    the red end of the strand along a straight line to the black vertex,
    thus turning the regular strand into a singular strand.
    Table~\ref{pinched:face_relation} shows a $\isstrictfaceof$-chain of
    length two.

    \begin{table}
      \begin{center}
        \setlength{\unitlength}{1mm}
        \begin{picture}(120,138.654)(0,-69.282)
          \drawline(0,0)(120,0)
          \drawline(60,34.641)(0,-69.282)
          \drawline(60,-34.641)(0,69.282)
          \drawline(0,-69.282)(0,69.282)(120,0)(0,-69.282)
          \put(40,0){\centergraphics[scale=0.5]{pinched_braid.eps}{}{}}
          \put(60,34.641){\centergraphics[scale=0.5]{pinched_face_one.eps}{}{}}
          \put(60,-34.641){\centergraphics[scale=0.5]{pinched_face_two.eps}{}{}}
          \put(0,0){\centergraphics[scale=0.5]{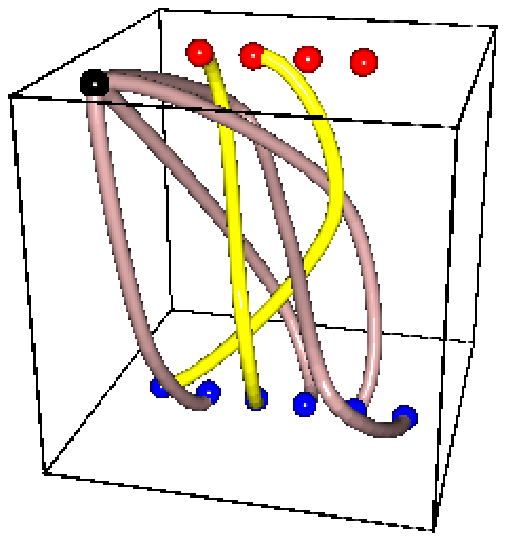}{}{}}
          \put(0,-69.282){\centergraphics[scale=0.5]{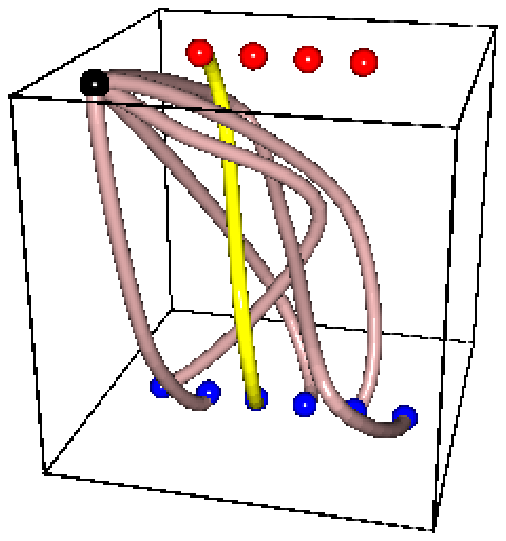}{}{}}
          \put(0,69.282){\centergraphics[scale=0.5]{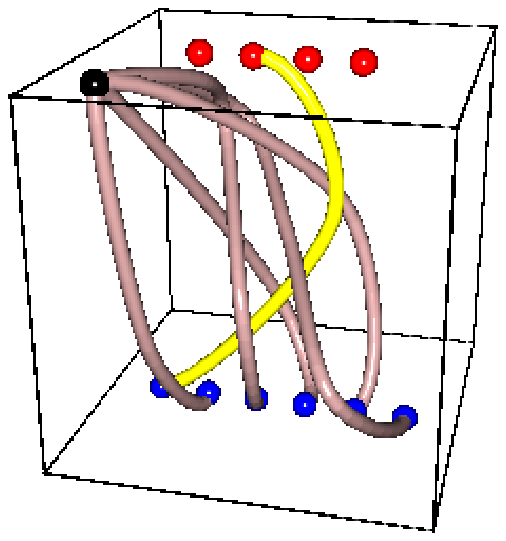}{}{}}
          \put(120,0){\centergraphics[scale=0.5]{pinched_vertex_three.eps}{}{}}
        \end{picture}
      \end{center}
      \bigskip
      \null
      \caption{\label{pinched:triangle}A triangle
        in $\PinchedBraidComplex[\Six][\Four]$ and its geometric realization}
    \end{table}
    The poset $\PinchedBraidPoset[\TheBotNumber][\TheTopNumber]$ thus defined
    is the poset of cells in a $\Delta$-complex whose vertex set is the set
    of pinched braids that have precisely one regular strand.
    Table~\ref{pinched:triangle} shows a $\Two$-simplex in this complex
    with a complete labelling of all its faces by pinched braids.
 
    We want to provide some evidence for the following:
    \begin{conj}\label{pinched:conj}
      $\PinchedBraidPoset[\TheBotNumber][\TheTopNumber]$
      is $\Parentheses{\TheTopNumber-\One}$-spherical
      provided $\TheBotNumber$ is large enough.
    \end{conj}
    \begin{rem}
      Deleting all singular strands defines a hight-preserving poset map
      \[
        \PinchedBraidPoset[\TheBotNumber][\TheTopNumber]
        \rightarrow
        \TheBraidsPoset[\TheBotNumber][\TheTopNumber]
        .
      \]
      Thus, a natural idea is to use Quillen's result. This
      thought led me to consider the braided chessboard complex
      in the first place. Unfortunately,
      the fibers of this projection seem to be difficult to
      analyze.
    \end{rem}
    \begin{rem}
      F.~Degenhardt
      \cite{Degenhardt:2000}
      introduced the series $\HoughtonBraid[\TheTopNumber]$ of
      \notion{braided Houghton groups}, and proved:
      \begin{enumerate}
        \item
          $\HoughtonBraid[\One]$ is not finitely generated.
        \item
          $\HoughtonBraid[\Two]$ is finitely generated but not
          finitely presented.
        \item
          $\HoughtonBraid[\Three]$ is finitely presented but not of
          type \FPType[\Three].
        \item
          $\HoughtonBraid[\TheTopNumber]$ is of type \FType[\Three] for
          $\TheTopNumber\geq\Four$.
      \end{enumerate}
      He conjectures that $\HoughtonBraid[\TheTopNumber]$ is of type
      \FType[\TheTopNumber-\One] but not of type \FType[\TheTopNumber].
      In an attempt to prove his conjecture, I constructed a contractible
      cube complex upon which $\HoughtonBraid[\TheTopNumber]$ acts
      with cell stabilizers of type $\FType[\infty]$.
      The complexes $\PinchedBraidPoset[\TheBotNumber][\TheTopNumber]$
      occur as relative links in a cocompact filtration by invariant
      subspaces. Thus, by standard arguments,
      Conjecture~\ref{pinched:conj} implies Degenhardt's conjecture
      on the finiteness properties of $\HoughtonBraid[\TheTopNumber]$.

      We remark that the Houghton groups $\Houghton[\TheTopNumber]$
      are groups of certain infinite
      permutations, and braided Houghton groups are groups of certain
      infinite braids. Ignoring the braiding defines a group homomorphism
      \(
        \HoughtonBraid[\TheTopNumber]
        \epirightarrow
        \Houghton[\TheTopNumber]
        .
      \)
      K.~Brown \cite[Section~5]{Brown:1989} derived the finiteness properties of
      $\Houghton[\TheTopNumber]$ from a filtration where the chessboard
      complexes $\ChessboardComplex[\TheBotNumber][\TheTopNumber]$ arose
      as relative links.
    \end{rem}

    We will show that Conjecture~\ref{pinched:conj} holds ``in the limit'':
    Adding an unused red vertex $\TheTopNumber+\One$ to
    the right of the top row induces an inclusion
    \[
      \PinchedBraidPoset[\TheBotNumber][\TheTopNumber]
      \subset
      \PinchedBraidPoset[\TheBotNumber][\TheTopNumber+\One]
      .
    \]
    Adding a blue vertex $\TheBotNumber+\One$ to the left of the bottom row,
    we define an embedding
    \[
      \PinchedBraidPoset[\TheBotNumber][\TheTopNumber]
      \subset
      \PinchedBraidPoset[\TheBotNumber+\One][\TheTopNumber]
    \]
    as follows: We fix a path in the boundary of the cube from the black
    vertex to the new blue vertex. For any pinched braid in
    $\PinchedBraidPoset[\TheBotNumber][\TheTopNumber]$ we define its image
    by pushing the boundary path into the cube, thereby creating a singular
    strand to the new blue vertex. This process is compatible with the face
    relation in $\PinchedBraidPoset[\TheBotNumber][\TheTopNumber]$ and,
    therefore, defines a poset morphism.

    Put
    \[
      \PinchedBraidPoset[\infty][\TheTopNumber] :=
      \Union[\TheBotNumber=\One][\infty]{
        \PinchedBraidPoset[\TheBotNumber][\TheTopNumber]
      }
    \]
    and
    \[
      \PinchedBraidPoset[\infty][\infty]
      :=
      \Union[\TheTopNumber=\One][\infty]{
        \PinchedBraidPoset[\infty][\TheTopNumber]
      }
      =
      \Union[\TheTopNumber,\TheBotNumber\in\NNN]{
        \PinchedBraidPoset[\TheBotNumber][\TheTopNumber]
      }
      .
    \]
    \begin{theorem}\label{pinched:main}
      $\PinchedBraidPoset[\infty][\TheTopNumber]$ is
      $\Parentheses{\TheTopNumber-\One}$-spherical.
    \end{theorem}
    
    A vertex $\TheVertex\in\PinchedBraidPoset[\infty][\infty]$
    involves precisely one regular strand.
    Let
    \[
      \TheTopSlotOf{\TheVertex}
    \]
    be the label of its top slot and let
    \[
      \TheBotSlotOf{\TheVertex}
    \]
    be the label of its bottom slot. Extending affinely
    to simplices, we define two \notion{height functions}
    \[
      \TheTopSlot \mapcolon
      \PinchedBraidPoset[\infty][\infty] \rightarrow \NNN
    \]
    and
    \[
      \TheBotSlot \mapcolon
      \PinchedBraidPoset[\infty][\infty] \rightarrow \NNN
      .
    \]
    Since there are no horizontal edges,
    these height functions are Morse functions as
    defined in \cite{Bestvina.Brady:1997}.
    Note that $\PinchedBraidPoset[\infty][\TheTopNumber]$
    is the sublevel set
    \[
      \SetOf[{
        \ThePoint\in\PinchedBraidPoset[\infty][\infty]
      }]{
        \TheTopSlotOf{\ThePoint} \leq \TheTopNumber
      }
      .
    \]
    \begin{observation}
      Consider a sphere in
      \(
        \PinchedBraidPoset[\TheBotNumber][\TheTopNumber]
        \subset
        \PinchedBraidPoset[\TheBotNumber+\One][\TheTopNumber+\One]
        .
      \)
      In all its simplices, we can slide the top end of the singular
      strand based at the bottom vertex $\TheBotNumber+\One$ to the top
      slot $\TheTopNumber+\One$. The regular strand thus created
      serves defines the same vertex in all simplices of the given
      sphere and, therefore, serves as a cone point from which
      the whole sphere can be contracted. Thus, the inclusion
      \[
        \PinchedBraidPoset[\TheBotNumber][\TheTopNumber]
        \monorightarrow
        \PinchedBraidPoset[\TheBotNumber+\One][\TheTopNumber+\One]
      \]
      is trivial in homotopy.\qed
    \end{observation}
    \begin{observation}
      Since any sphere in $\PinchedBraidPoset[\infty][\TheTopNumber]$
      involves only finitely many
      cells, the argument just given also implies that the inclusion
      \[
        \PinchedBraidPoset[\infty][\TheTopNumber]
        \monorightarrow
        \PinchedBraidPoset[\infty][\TheTopNumber+\One]
      \]
      is trivial in homotopy.

      In particular, all homotopy groups of
      $\PinchedBraidPoset[\infty][\infty]$ vanish, i.e.,
      $\PinchedBraidPoset[\infty][\infty]$ is
      contractible.\qed
    \end{observation}

    We need to generalize Theorem~\ref{pinched:main} a little to
    make it amenable to an inductive argument.
    Let $\PinchedBraidPosetOf[\TheBotNumber][\TheTopNumber]{\TheBarNumber}$
    be the poset of pinched
    braids with $\TheBotNumber$ blue bottom vertices, $\TheTopNumber$
    red top vertices   , one $\infty$-slot in front of the
    top row, and $\TheBarNumber$ green fixed disjoint vertical strands
    connecting $\TheBarNumber$ \emph{additional} pairs of
    vertices. These green strands are not involved
    in the definition of the face relation, they stay
    put. This generalization now describes a class of
    complexes closed with respect to taking links:
    \begin{observation}
      The link of a vertex in
      \(
        \PinchedBraidPosetOf[\TheBotNumber+\One][\TheTopNumber+\One]{
          \TheBarNumber
        }
      \)
      is
      isomorphic to
      \(
        \PinchedBraidPosetOf[\TheBotNumber][\TheTopNumber]{
          \TheBarNumber+\One
        }
        .
      \)\qed
    \end{observation}

    Note that our previous observation carries over to the more general setting:
    \begin{observation}\label{pinched:obs}
      The inclusion
      \[
        \PinchedBraidPosetOf[\infty][\TheTopNumber]{\TheBarNumber}
        \monorightarrow
        \PinchedBraidPosetOf[\infty][\TheTopNumber+\One]{\TheBarNumber}
      \]
      is trivial in homotopy.\qed
    \end{observation}

    The following includes Theorem~\ref{pinched:main}:
    \begin{theorem}
      The map
      \[
        \HomotopyOf[\TheIndex]{
          \PinchedBraidPosetOf[\infty][\TheTopNumber]{\TheBarNumber}
        }
        \rightarrow
        \HomotopyOf[\TheIndex]{
          \PinchedBraidPosetOf[\infty][\TheTopNumber+\One]{\TheBarNumber}
        }
      \]
      induced by the inclusion is an isomorphism
      for $\TheIndex<\TheTopNumber-\Two$ and onto
      for $\TheIndex=\TheTopNumber-\Two$. In particular, the space
      \(
        \PinchedBraidPosetOf[\infty][\TheTopNumber]{
          \TheBarNumber
        }
      \)
      is $\Parentheses{\TheTopNumber-\Two}$-connected in view of
      Observation~\ref{pinched:obs}.
    \end{theorem}
    \begin{proof}
      This is combinatorial Morse theory and
      induction on $\TheTopNumber$: Consider the height function
      \[
        \TheBotSlot \mapcolon
        \PinchedBraidPosetOf[\infty][\infty]{
          \TheBarNumber
        }
        \rightarrow \NNN
        .
      \]
      The descending links of vertices of
      height $\TheTopNumber$ are isomorphic to
      \(
        \PinchedBraidPosetOf[\infty][\TheTopNumber-\One]{
          \TheBarNumber+\One
        }
        .
      \)
      This complex is $\Parentheses{\TheTopNumber-\Three}$-connected
       by induction. Thus the
      statement follows
      from \cite[Lemma~2.5 and Corollary.6]{Bestvina.Brady:1997}.
    \end{proof}

  \section{References}
  \InputIfFileExists{bibliography.tex}{}{}
\end{document}